\documentclass[]{interact}

\usepackage{epstopdf}
\usepackage[caption=false]{subfig}

\usepackage[numbers,sort&compress]{natbib}
\bibpunct[, ]{[}{]}{,}{n}{,}{,}
\renewcommand\bibfont{\fontsize{10}{12}\selectfont}
\makeatletter
\def\NAT@def@citea{\def\@citea{\NAT@separator}}
\makeatother

\theoremstyle{plain}

\theoremstyle{definition}

\newtheorem{thm}{Theorem}
\newtheorem{prp}{Proposition}

\newtheorem{cor}{Corollary}
\newtheorem{rem}{Remark}

\theoremstyle{remark}

\begin{document}

\title{Complex invariants of poristic Steiner $4$-chains}

\author{
	\name{Ana Diakvnishvili\textsuperscript{a}\thanks{CONTACT Ana Diakvnishvili. Email: ana.diakvnishvili.1@iliauni.edu.ge, Giorgi Khimshiashvili. Email: giorgi.khimshiashvili@iliauni.edu.ge} and Giorgi Khimshiashvili\textsuperscript{b}}
	\affil{\textsuperscript{a,b}Faculty of Business,Technology and Education, Ilia State University, Tbilisi, Georgia}}

\maketitle

\begin{abstract}
	
	We are concerned with the Steiner chains consisting of four circles. More precisely, we deal with the so-called complex moments of Steiner $4$-chains introduced in a recent paper by J.Lagarias, C.Mallows and A.Wilks. We compute the invariant complex moments of poristic Steiner $4$-chains and establish certain algebraic relations between those invariants. To this end we use the invariance of certain moments of curvatures of poristic Steiner chains established by R.Schwartz and S.Tabachnikov, combined with the computation of these moments for the so-called symmetric Steiner $4$-chains. We also present analogous results for poristic Steiner $3$-chains and give an application to the feasibility problem for the centers of Steiner $4$-chains.
\end{abstract}


\vspace*{0.25cm}\noindent{\small {\bf MSC 2010:} {52C35,  32S40.}}

\section{Introduction}\label{sec1}

We present a number of results on the complex moments of poristic Steiner $4$-chains introduced in \cite{lamawi} (cf. also \cite{scta}), based on the results on the so-called symmetric Steiner chains given in \cite{bibdia}. In particular, we compute the invariant complex moments through the radii of the parent circles of poristic Steiner $4$-chains and give certain algebraic relations between those invariants. Our approach is basically the same as in \cite{bibdia} and relies on a few general results on the Steiner chains given in \cite{scta} and \cite{yiu}. \\

In the second section we present the necessary background on the invariants of Steiner chains introduced in \cite{lamawi} and \cite{scta}. In the third section we recall a number of general results about symmetric Steiner chains and compute some of the invariants defined in the second section. In the fourth section we establish certain algebraic relations between the invariants considered and outline an application of the obtained results to the feasibility problem for Steiner $4$-chains. We also present a few similar results for Steiner $3$-chains. The last section contains several remarks on possible generalizations and research perspectives suggested by our results.\\


We proceed by recalling that the term {\it Steiner $n$-chain} (aka $n$-garland) refers to a sequence of circles $\delta_i, i=1, \ldots, n,$ in Euclidean plane such that the circles with adjacent indices $mod\, n$ are externally tangent, each of $\delta_i$ is externally tangent to a fixed circle $\gamma$ and internally tangent to another fixed circle $\Gamma$ containing $\gamma$ in its interior disk. The nested circles $\gamma$ and $\Gamma$ are called the {\it parent circles} of the Steiner chain considered \cite{ber}. Steiner chains is a classical topic discussed in many papers, in particular, in the context of {\it Steiner porism} which is in the focus of our discussion (see, e.g., \cite{ped,scta,bibdia,dia},). \\

We note that triples of pairwise externally tangent circles are often called "kissing circles". An ancient result of Apollonius states that, for any triple of "kissing circles", there exist exactly two circles which are tangent to each of the three given circles. These two circles are called the {\it parent circles} of the given kissing circles. The famous Descartes' circle theorem enables one to compute the curvatures of parent circles in terms of the curvatures of the three given "kissing circles" \cite{scta}. \\

The aforementioned Steiner porism states that there exists a one-dimensional family of Steiner $n$-chains having the same parent circles \cite{scta}. By a way of analogy with the Poncelet porism \cite{ber} the collection of Steiner $3$-chains with the fixed parent circles will be called {\it poristic Steiner $n$-chains}. This setting suggests, in particular, several extremal problems for poristic Steiner $n$-chains with the fixed Soddy circles. In particular, a natural extremal problem concerned with the sum of areas of poristic Steiner $4$-chains was studied in a prize-winning paper of K.Kiradjiev \cite{kir}. The results of K.Kirajiev have been generalized in a joint paper of G.Bibileishvili with the author of the present paper \cite{bibdia}. \\

The approach and results of \cite{bibdia} served as a motivation and impetus for the research which led to the results described in the present paper.

	\section{Moments and invariants of Steiner chains} \label{invariantmoments}

Let $n\geq 3$ be a natural number and $(R, r, d)$ with $R>r$ be a triple of positive numbers satisfying the relation 	
\begin{equation} \label{pedoe-n}
	d^2 =(R-r)^2-4qRr,
\end{equation}
where $q=\tan^2\frac{\pi}{n}.$ It is known that such a triple defines a poristic system of Steiner $n$-chains (aka Steiner $n$-garlands) with the parent circles having radii $R, r$ and the distance between their centers equal to $d$ (see, e.g., \cite{ped}). We will say that this is a pair of parent circles of order $n$ with the gauge $(R,r,d)$. The number $d$ will be called the {\it gap} of parent circles. Denote by $a= \frac{1}{r}, A=-\frac{1}{R}$ the so-called {\it bends} (signed curvatures) of the parent circles. Following \cite{scta}, for each Steiner $n$-chain $G$ and natural number $k$, we define the $k$-th {\it bending moment} $I_k^{(n)}(G)$ of chain $G$ as
\begin{equation} \label{bendingmoments1}
	I_k^{(n)}(G) = \sum_{j=1}^n b_j^k,
\end{equation}
where $b_j$ are the bends (curvatures) of the circles in $G.$ As was shown in \cite{scta}, the first $n-1$ bending moments $I_k^{(n)}(G)$ are invariant in a poristic family of poristic Steiner $n$-chains with the given parent circles, so we denote them simply by $I_k^{(n)} = I_k^{(n)}(R,r,d)$ and call $I_k^{(n)}$ the {\it invariant bending moments} of poristic Steiner chains. Clearly, they can be expressed through radii $R, r$ of the parent circles or, equivalently, through their bends (signed curvatures) $a, A$. We will give such formulas for $n=3$ and $n=4$. \\

If the gap $d$ is equal to the radius $r$ of the inner parent circle of order $n$, such a pair of parent circles  will be called {\it exact}. Clearly, this is equivalent to the relation $R = (2q+4)r$ with $q=\tan^2\frac{\pi}{n}.$ So the invariants of the corresponding poristic Steiner chains only depend on $r$, which enables one to obtain quite detailed results in this case and gives important hints for the general case, where $d\neq r$. For this reason the exact parent circles serve for us as running examples and play an important role in the sequel.\\

For $n=3$, explicit formulas for the first two bending moments $I_1, I_2$ follow directly from the Descartes circle theorem (see, e.g., \cite{scta}):
\begin{equation} \label{moment31}
	I_1^{(3)} = b_1 + b_2 + b_3 = \frac{A+a}{2} = \frac{R-r}{2Rr},
\end{equation}
\begin{equation} \label{moment32}
	I_2^{(3)} = b_1^2 + b_2^2 + b_3^2 = \frac{A^2 + 6Aa + a^2}{8} = \\
	\frac{R^2 - 6Rr + r^2}{8R^2r^2}.\\		.
\end{equation}

\begin{rem}
	It should be noted that the formula for $I_2^{(3)}$ given in page 239 of \cite{scta} is not correct and should be replaced by (\ref{moment32}) which has been proved in \cite{bibdia}.
\end{rem}

For $n=4$, analogous explicit formulas for the invariant bending moments $I_j^{(4)}, j=1,2,3,$ have been obtained in \cite{bibdia} using a general result given in \cite{yiu}. It should be noted that this task was not considered neither in \cite{scta}, nor in \cite{yiu}. According to \cite{bibdia} we have:

\begin{equation} \label{moment41}
	I_1^{(4)} = \frac{2(R-r)}{Rr} = 2(A+a),
\end{equation}

\begin{equation} \label{moment42}
	I_2^{(4)} = \frac{3R^2 - 10Rr + 3r^2}{2R^2r^2} = \frac{3A^2 + 10Aa + 3a^2}{2},
\end{equation}

\begin{equation} \label{moment43}
	I_3^{(4)} = \frac{5R^3 - 27R^2r + 27Rr^2 - 5r^3}{4R^3r^3} = \frac{5A^3 + 27A^2a + 27Aa^2 + 5a^3}{4}.
\end{equation}

We recall now the aforementioned general result on Steiner $n$-chains established in \cite{yiu} since it will be used in the sequel. 

\begin{prp} \label{eqYiu} (\cite{yiu})
	Let $C(u)$ be a poristic circle in a Steiner $n$-chain having radius $u$ and let $v = \frac{1}{u}$ denote its curvature. Then the curvatures 
	$v_+, v_-$ of the two neighbours of $C(u)$ are the roots of quadratic trinomial
	\begin{equation} \label{yiuquadratics}
		\alpha x^2 + \beta x + \gamma,
	\end{equation}
	where
	\begin{equation} \label{yiucoefficients}
		\begin{split}
			\alpha = (q+1)^2R^2r^2u^2, \\
			\beta = 2(q+1)Rru[(q-1)Rr - (R-r)u], \\
			\gamma = [(q+1)Rr - (R-r)u]^2 + 4Rru^2. \\
		\end{split}
	\end{equation}
\end{prp}

The following corollary is obtained by merely applying Bezout theorem to the above equation.

\begin{cor}
	The sum $v_+ + v_-$ of two neighboring bends of $C(u)$ is equal to $-\frac{\beta}{\alpha}.$
\end{cor}

Changing the variable $x$ by its inverse we get an analogous conclusion for the neighboring radii which will be used in the sequel.

\begin{cor}
	The sum $u_+ + u_-$ of two neighboring radii of $C(u)$ is equal to $-\frac{\beta}{\gamma}.$
\end{cor}

We will also use the following general property of Steiner $n$-chains presented in \cite{bibdia}.

\begin{prp} \label{poristicrange1}
	For a given pair of parent circles with the gauge $(R, r, d)$, the minimal and maximal possible values of poristic circles $r_i$ are
	$$r_* = \frac{R - d - r}{2}, r^* = \frac{R + d - r}{2},$$
	while the minimal and maximal values of poristic curvatures are
	$$b_* = \frac{2}{R-d-r}, b^* = \frac{2}{R+d-r}$$
	respectively. For any $r \in [r_*, r^*],$ the poristic family contains a circle of radius $r$. 
\end{prp}

\begin{rem}
	We notice that a Steiner chain with the center of one its circles belonging to $L$ is invariant with respect to the reflection in $L$. Chains with this property will be called {\it symmetric Steiner chains}. Symmetric Steiner chains will play an important role in the sequel.
\end{rem}  

Proposition \ref{eqYiu} enables one to compute the radii and bends of the circles in a Steiner garland. In the sequel we also need to compute the centers of poristic circles. To this end we present another general result on Steiner chains, which is possibly new since we could not find it in the existing literature. \\

Recall that the term {\it axis of porism} refers to the straight line through the centers of parent circles. Obviously, this line is well defined only in non-concentric case, i.e. if $d\neq 0$. In concentric case any line though the common center of parent circles can be regarded as an axis of porism. Let us introduce a coordinate system with the origin at the center of the outer parent circle and $Ox$ axis along the axis of porism. Such a coordinate system will be called the {\it canonical coordinate system} of a porism. We give an exact formula for the radius $r(t)$ of poristic circle with the center $z(t)$ having a given polar angle $t$ in the canonical coordinate system.  

\begin{prp} \label{poristicradius}
	The radius $r(t)$ of a poristic circle with the center having a given polar angle $t$ in the canonical coordinate system is given by:	
	\begin{equation} \label{keyformula1}	
		r(t) = \frac{R^2 - 2dR\cos t + d^2 - r^2}{2(R + r - d\cos t)}.
	\end{equation}
\end{prp}

{\bf Proof.} We apply the cosine rule to the triangle formed by the centers $\Omega$ and $\omega$ of parent circles and the center $z(t)$ of the poristic triangle considered. From the very definition of Steiner chain follows that the sides of this triangle are $R - r(t)$, $r + r(t)$ and $d$. The result follows by noticing that the side $\vert\omega z(t)\vert$ is opposite to the angle $t$ and resolving the cosine rule with respect to $r(t)$. \\

This proposition will be repeatedly used in the sequel as well as the following corollary which follows by inverting equation (\ref{keyformula1}).

\begin{cor} \label{polarangle}
	For any $\rho \in [r_*, r^*],$ there exist exactly two poristic circles with the radii equal to $\rho$ and their polar angles $t_{\pm}(\rho)$ are given by  
	\begin{equation} \label{keyformula2}	
		t_{\pm}(\rho) = \pm \cos^{-1} \Big(\frac{R^2 - r^2 + d^2 - 2(R+r)\rho}{2d(R-r)}\Big).
	\end{equation}
\end{cor}

Another corollary yields the interval where $r(t)$ is unimodal, which implies the log-concavity property playing important role in various combinatorial issues. 

\begin{cor} \label{unimodality}
	The poristic radius $r(t)$ is unimodal in segments $[0, \pi]$ and $[\pi, 2\pi].$
\end{cor}	
This corollary is proved by noticing that the derivative $dr/dt$ of the right-hand-side of (\ref{keyformula1})  
is equal to $\sin t$ times a certain constant, which obviously implies the statement. \\

Using Proposition \ref{poristicradius} and Corollary \ref{polarangle} one can easily compute the (canonical) coordinates of the center $z(t)$ as follows. For simplicity, we only describe this procedure for $n=3$. As was already mentioned, the distances between $z(t)$ and the centers of the outer and inner parent circles are $R-r(t)$ and $r + r(t),$ respectively. Hence the abscissa and ordinate of $z(t)$ are equal to $(R - r(t)) \cos t$ and $(R - r(t)) \sin t$ respectively. The radii $r_{\pm}$ of two neighboring poristic circles can be found as the roots of quadratic polynomial (\ref{yiuquadratics}). Having the values of $r_{\pm}$, the polar angles of the centers $z_{\pm}$ of the two neighboring poristic circles are given by equation (\ref{keyformula2}). For $n=3$, this yields explicit formulas for the centers and radii of all poristic circles as functions of $t,$ which will be used in the sequel. For $n\geq 4,$ one can do the same inductively but the resulting formulas are more difficult 
to obtain. We will only give more details in the case where $n=4$ which is our main concern in this paper. \\

Following a suggestion given in \cite{lamawi}, for any natural $n\geq 3$, we consider further numerical characteristics $J_{k,m}(G)$ of a Steiner $n$-chain $G$ defined as follows. For a pair of non-negative integers $k, m$ and a Steiner $n$-chain $Z$, we define 
\begin{equation} \label{complexmoments}
	J_{k,m}(G) = \sum b_i^k\cdot z_i^m,
\end{equation}
where $b_i$ are the bends as above and $z_i$ are the centers of the circles in the given Steiner $n$-chain $G$, considered as complex numbers by complexifying the reference plane endowed with canonical coordinate system of porism. For brevity, we will refer to the quantities $J_{k,m}(G)$ as the {\it complex moments} of the Steiner chain $G$. Clearly, $J_{k,0} = I_k$ so the complex moments are generalizations of the bending moments. \\

An important result of \cite{lamawi} states that the complex moments $J_{k,m}$ with $0\leq m \leq k \leq n-1$ are invariants of poristic Steiner $n$-chains. We call such quantities the {\it invariant complex moments} of poristic Steiner chains and establish several results about those invariants. Clearly, the total amount of all invariant moments is $(n-1)(n+2)/2$ including $n-1$ bending moments $J_{k,0} = I_k.$ \\

Clearly, the quantities $J_{k,m}(G)$ are in general complex numbers. However the invariant complex moments have in fact a hidden real structure. Namely, in the next section we will show that all quantities $J_{k,m}$ with $0\leq m \leq k \leq n-1$ become simultaneously real for an appropriate embedding of the parent circles in the complex plane. We emphasize that this property only refers to invariant complex moments. These observations may be new since to the best of our knowledge they were never mentioned in the literature. \\

An essential novelty of the present paper, as compared with \cite{bibdia} and \cite{lamawi}, is that we also consider the invariant complex moments $J_{k,m}$ and moment loops $L_{k,m}$. Another novelty is that we prove the existence of algebraic relations between the invariant moments. Surprisingly, the latter topic was not explored in \cite{lamawi} and \cite{scta}. We also explain how one can use such relations for solving the so-called {\it feasibility problem} for Steiner $4$-chains considered in \cite{kir} and \cite{bisa}. An important role in all these developments is played by the symmetric Steiner chains discussed in the next section. \\

We conclude this section by computing the invariant complex moments for a concrete pair of parent circles of order three. The first example corresponds to an exact pair of parent circles while the second one arise from a generic pair of parent circles. We compute the invariant complex moments $J_{1,1}, J_{2,1}, J_{2,2}.$ As was explained above, the moment loop $L_{3,1}$ can be plotted using Proposition \ref{poristicradius} and its corollary. \\

{\bf Example 1.} Consider a pair of circles with the gauge $(14, 1, 1).$ It is easy to see that these circles satisfy the relation (\ref{pedoe-n}) for $n=3$ so they are the parent circles of a Steiner porism of order three. In this example $r^* = 9$ and computing the bends of the maximal $3$-chain $C^*$ one gets $b_1 = 0.111, b_2 = b_3 = 0.1775.$ Introducing the system of coordinates with the origin at the center of outer parent circle and having the values of axial radii it is easy to see that the centers of the circles in $C^*$ are 
$z_1 = (-6,0), z_2 = (7.5, 7.66), z_3 = (7.5,-7.66)$. After that it is straightforward to compute the values of the invariants discussed above and we get that they are $6.5 ,4.3, 3.7$ respectively. \\

{\bf Example 2.} Consider a pair of circles with the gauge $(15, 1, 4).$ It is easy to see that these circles satisfy the relation (\ref{pedoe-n}) for $n=3$ so they are the parent circles of a Steiner porism of order three. In this example $r^* = 9$ and computing the bends of the maximal $3$-chain $C^*$ one gets $b_1 = 0.111, b_2 = b_3 = 0.1775.$ Introducing the system of coordinates with the origin at the center of outer parent circle and having the values of axial radii it is easy to see that the centers of the circles in $C^*$ are 
$z_1 = (-6,0), z_2 = (7.5, 7.66), z_3 = (7.5,-7.66)$. After that it is straightforward to compute the values of the invariants discussed above and we have 
$$I_1 = \frac{14}{30} = 0.466, I_2 = 0.074, J_{11} = 2.$$

Further examples of computations of invariant bending moments for $n=3$ can be found in \cite{bisa}. 
The invariant complex moments for Steiner $4$-chains will be computed in the sequel using the properties of symmetric Steiner chains discussed in the next section.

	\section{Symmetric Steiner chains } \label{symmetries}

For arbitrary $n\geq{3}$ and a pair of parent circles of Steiner porism of order $n$, poristic Steiner chains contain the symmetric chains mentioned above. In what follows it is convenient to introduce a Cartesian coordinate system in the reference plane such that $Ox$-axis coincides with the axis of porism and the origin coincides with the center of outer parent circle. This system of coordinates will be called the canonical coordinate system associated with the given pair of parent circles and the coordinates of a point in this system will be called its {\it canonical coordinates}. Consideration of such coordinates enables us to identify the reference plane with the complex plane $\mathbb{C}$. In this coordinate system, the symmetry of a Steiner chain is equivalent to its invariance with respect to the complex conjugation. It is also convenient to distinguish between the two types of symmetric configurations: axial and lateral ones. \\

It is easy to see that among poristic Steiner $n$-chains there always are a chain $C^*$ containing a circle with the biggest possible radius $r^*$ and another one, $C_*,$ containing a circle with the smallest possible radius $r_*$ given above. Both these chains are called axial $n$-chains since they are symmetric with respect to the axis of porism. For even $n,$ the two chains $C_*$ and $C^*$ coincide and will be denoted by $C_*^*$. For even $n$, there also exists a symmetric chain having only two common points with the axis at which it is tangent to the axis of porism. The latter symmetric chain will be called the {\it lateral symmetric chain} of poristic family. Thus, for odd $n$, there are always two different axial symmetric configurations, while for even $n$, there exist one axial symmetric chain  $C_*^*$ and one lateral symmetric chain. \\

For $n=3,$ the explicit formulas for the radii and bends of the axial chains have been given in \cite{bisa}. We complement the above result of \cite{bisa} by giving analogous formulas for $n=4,$ \cite{bibdia} which are proved using equation (\ref{yiuquadratics}). In this case $q=1$. For the biggest circle in $C_*^*$ we have again $$u = r^* = \frac{R-d+r}{2},v = \frac{2}{R-r+d}$$ and the coefficients of equation (\ref{yiuquadratics}) take the form:
\begin{equation} \label{yiu4}
	\alpha = 4R^2r^2u^2, \beta = -4Rru^2(R-r),
	\gamma = [2Rr - (R-r)u]^2 + 4Rru^2.
\end{equation}
Solving the latter equation we get the values of radii and bends in the axial chain. The Yiu's equation  (\ref{yiuquadratics}) also enables one to find the radii and bends in the lateral symmetric garland for $n=4$. Denote the two unknown lateral curvatures by $b_+$ and $b_-$. From the invariance of the first two bending moments and the formulas for $I_1^{(4)}, I_2^{(4)}$ given above we get a system of equations for $b_+, b_-:$
$$\{b_+ + b_- = A+a, b_+^2 + b_-^2 = \frac{3A^2 + 10Aa + 3a^2}{4}\}.$$
It follows that $b_+, b_-$ are the roots of quadratic equation
$$x^2 - (A+a)x + \frac{(A-a)^2}{8} = 0$$
solving which the find the radii in the lateral chain. In this way we arrive to the following result. 	

\begin{prp} 
	For a given pair of parent circles of order four with the gauge $(R,r,d)$, there are both the axial and lateral $4$-chains. The quadruples of their radii $r_i$ and bends $b_i$ are given below. \\
	The axial quadruple of radii is:
	\begin{equation} \label{4axialradii}
		\left(\frac{R - d - r}{2}, \frac{2Rr}{R-r}, \frac{R + d - r}{2}, \frac{2Rr}{R-r}\right).
	\end{equation} \\ 	 	
	Correspondingly, the axial quadruple of bends is:
	\begin{equation} \label{4axialbends}
		\left(\frac{2}{R - d - r}, \frac{R-r}{2Rr}, \frac{2}{R + d - r}, \frac{R-r}{2Rr}\right).
	\end{equation} \\
	The lateral quadruple of bends is: 
	\begin{equation} \label{4lateralbends}
		\left(\frac{R-r}{2Rr} - \frac{d}{2Rr\sqrt{2}}, \frac{R-r}{2Rr} + \frac{d}{2Rr\sqrt{2}},\\ 
		\frac{R-r}{2Rr} + \frac{d}{2Rr\sqrt{2}}, \frac{R-r}{2Rr} - \frac{d}{2Rr\sqrt{2}}\right).					
	\end{equation}
\end{prp}
The lateral quadruple of radii is: 
\begin{equation} \label{4lateralradii}
	\left(\dfrac{4 R r}{2(R - r) - \sqrt{2}\, d},
		\dfrac{4 R r}{2(R - r) + \sqrt{2}\, d}, \\
		\dfrac{4 R r}{2(R - r) + \sqrt{2}\, d},
		\dfrac{4 R r}{2(R - r) - \sqrt{2}\, d}\right)					
\end{equation}

Now, using the procedure described in Section 2 one can explicitly compute the centers of circles in the axial Steiner $4$-chain, which enables us to compute the invariant complex moments $J_{k,m}$ using formula (\ref{complexmoments}). In this way, we get a simple expression for the invariant complex moments:

\begin{thm}
	For a given pair of parent circles of order four with the gauge $(R,r,d)$, we have
\begin{equation}
	J^{(4)}_{1,1} = \frac{2d}{r}=\frac{2\sqrt{A^{2}+6Aa+a^{2}}}{A}
\end{equation}

\begin{equation}
	J^{(4)}_{2,1} = \frac{d\,(3R - 5r)}{2 R r}=\frac{3a+5A}{2A}\,\sqrt{A^{2}+6Aa+a^{2}}
\end{equation}

\begin{equation}
	J^{(4)}_{3,1} = \frac{(5R - 3r)(R - 3r)\,d}{4 R^{2} r^{3}}=\frac{(5a+3A)(a+3A)}{4A}\,\sqrt{A^{2}+6Aa+a^{2}}
\end{equation}

\begin{equation}
	J^{(4)}_{2,2} = \frac{3 d^{2}}{2r}=\frac{3\sqrt{A^{2}+6Aa+a^{2}}}{4A}
\end{equation}
\end{thm}  

It turns out that properties of symmetric Steiner chains enable us to get another general result about the invariant complex moments. To this end consider the canonical coordinate system introduced above. Then the invariance of axial chains under the complex conjugation implies that the values of all the invariant complex moments computed at an axial chain are real numbers. This yields the following result which reveals the hidden real structure of the invariant complex moments.  

\begin{prp} \label{invariantcomplexmoments2}
	The invariant complex moments are real in the canonical coordinate system.	
\end{prp}

It is easy to see that this conclusion does not refer to the higher complex moments. In other words,
higher complex moments are usually complex numbers, which can be seen in simple examples. \\

In the next section we show that there exist algebraic relations between the invariants introduced above and outline their application to the so-called {\it feasibility problem} for Steiner $4$-chains.  

\section{\textbf{Algebraic relations between invariants}}

We proceed by discussing the algebraic relations between the invariant moments. As will be explained in the proof of Theorem \ref{main2} below, the existence of such relations follows from some general results. To make our exposition more clear we begin by discussing in some detail the first of such relations for $n=4$. According to \cite{scta}, in this case there are three invariant bending moments $I_i, i=1,2,3$ and they all are expressed by symmetric polynomials in two variables $(A, a)$. By a fundamental algebraic theorem on algebraic dependence of polynomials there exists an algebraic relation between $I_1, I_2$ and $I_3$. It turns out that in our case this relation is especially simple. Namely, the third bending moment is expressed as a polynomial in the first two bending moments. The exact form of this relation was given in \cite{dia}.
\begin{thm} \label{thirdmomentrelation}
	With the above notations and assumptions one has:
	\begin{equation} \label{moment3byfirst2}
		I_{3}^{(4)}=\frac{3}{4}I_1I_2-\frac{1}{8}I_{1} ^{3}. 
	\end{equation}
\end{thm}
{\bf Proof.} From the formulas for $I_1^{(4)}$, $I_2^{(4)}$ and $I_3^{(4)}$ given above we have:

$$I_1^{(4)} = 2(A + a),  A + a = \frac{I_1^{(4)}}{2}.$$
We also have
$$I_2^{(4)} = \frac{3A^2 + 10Aa + 3a^2}{2}
= \frac{3(A + a)^2 + 4Aa}{2}
= \frac{3\left(\frac{I_1^{(4)}}{2}\right)^2 + 4Aa}{2}.$$
Hence  
$$I_2^{(4)} = \frac{3}{8}(I_1^{(4)})^2 + 2Aa,$$
from where we express $Aa$ as

$$Aa = \frac{I_2^{(4)}}{2} - \frac{3}{16}(I_1^{(4)})^2.$$
We can now rewrite the expression for the third moment as follows;

$$I_3^{(4)} = \frac{5A^3 + 27A^2a + 27Aa^2 + 5a^3}{4}$$

$$= \frac{5A^3 + 15A^2a + 15Aa^2 + 12A^2a + 12Aa^2 + 5a^3}{4}$$

$$= \frac{5(A + a)^3 + 12Aa(A + a)}{4}.$$
Finally, inserting the expressions $A+a$ and $Aa$ into the latter formula we get

$$I_3^{(4)} = \frac{5\left( \frac{I_1^{(4)}}{2} \right)^3 + 12\left( \frac{I_2^{(4)}}{2} - \frac{3}{16}(I_1^{(4)})^2 \right)\left( \frac{I_1^{(4)}}{2} \right)}{4}$$

$$= \frac{3}{4}I_1^{(4)}I_2^{(4)} - \frac{1}{8}(I_1^{(4)})^3,$$
which completes the proof of the theorem.\\

Similar relations are available for complex invariant moments. For $J_{3,1}^{(4)},$ it immediately follows for the formulas (18),(19),(20) by squaring and using the Pedoe relation for $n=4$.   		
\begin{prp} \label{complexmomentrelation1}
	With the squaring and using the Pedoe relation we have:
	\begin{equation}
		J_{1,1}^2 =\frac{4\left(A^{2}+6Aa+a^{2}\right)}{A^{2}} ,
	\end{equation}
    \begin{equation}
	    J_{2,1}^2 =\frac{(3a+5A)^{2}\left(A^{2}+6Aa+a^{2}\right)}{4A^{2}},
    \end{equation}
    \begin{equation}
	    J_{3,1}^2 =\frac{(5a+3A)^{2}(a+3A)^{2}\left(A^{2}+6Aa+a^{2}\right)}{16A^{2}},
    \end{equation}
    \begin{equation}
	J_{2,2}^2 =\frac{9\left(A^{2}+6Aa+a^{2}\right)^{2}}{4A^{4}}.
    \end{equation}
\end{prp}

Analogous relations for other invariant complex moments can be obtained in similar way by iterated squaring but they are much more complicated and therefore omitted.

\begin{thm}
	With the above notations and assumptions one has:
	\begin{equation}
	J_{3,1}^{(4)}=-\frac{4\left(-J_{1,1}^{2}+12M-68\right)^{3} J_{2,1}^{4}
		\left(-5J_{1,1}^{2}+6M-16\right)^{2}
		\left(-J_{1,1}^{2}+6M-32\right)^{2}}{J_{1,1}^{2}
		K^{2}
		(M-6)^{2}}	
	\end{equation}
where $M=\sqrt{J_{1,1}^{2}+32},$\\
$K=\left(-9J_{1,1}^{6}-1720J_{1,1}^{4}+168M J_{1,1}^{4}-66448J_{1,1}^{2}+9696M J_{1,1}^{2}+131328M-742912\right)$
\end{thm}

\begin{rem}
	As will be shown below formula (\ref{moment3byfirst2}), in particular, can be used to elaborate upon a criterion for virtual quadruples of radii of Steiner $4$-chains given in \cite{bibdia}. It also suggests further applications of this kind involving higher bending moments.
\end{rem}

The general result on the existence of algebraic relations between invariant bending and complex moments follows by a similar argument using an explicit formulas for the first two bending moments communicated to us by M.Jibladze. For the convenience of further reference, these formulas are given in a separate proposition. 

\begin{thm} \label{first2moments}
	For any $n\geq 3,$ with the above notations and assumptions the following formulas hold true:
	\begin{equation} \label{mamuka1}
		I_1^{(n)} = ns, I_2^{(n)} = \frac{n}{2}(3s^2 + p),
	\end{equation}
	where 
	$$s = cot^2(\frac{\pi}{n})\frac{A+a}{2}, p = cot^2(\frac{\pi}{n})Aa.$$ 
\end{thm}

From the formulas (\ref{mamuka1}) follows that, for any $n\geq 3$, one can express $a+A$ and $aA$ as rational functions of $I_1$ and $I_2.$ Since all invariant bending moments are symmetric polynomials in $A$ and $a$, this implies that all invariant bending moments, which by a fundamental theorem on symmetric polynomials are polynomials in $a+A$ and $aA$, can be algebraically expressed through $I_1$ and $I_2.$ \\

Since the bends of axial configurations are also algebraically expressible through the bends of parent circles, formulas (\ref{mamuka1}) imply that the complex invariant moments are also algebraically expressible through the parent bends $A$ and $a$. So we have the following general conclusion. 
\begin{thm} \label{main2}
	For any $n\geq 4$, the invariant bending moments $I_3, \ldots, I_{n-1}$ of poristic Steiner $n$-chains are algebraically expressible through the first two bending moments $I_1, I_2$. 
\end{thm}

\begin{rem}
	For any concrete $n$, one can obtain the explicit form of the aforementioned relations by using the standard results on symmetric polynomials. However it seems rather difficult to obtain general formulas of such kind. As M.Jibladze informed us, he obtained analogous general formulas for all invariant bending moments using an inversion which transforms the given parent circles into a pair of concentric ones. It would be interesting to investigate if the formulas mentioned by M.Jibladze may be used to obtain explicit relations between invariant complex moments for arbitrary $n$. We emphasize that our results for Steiner $4$-chains are proved independently since, for $n\leq 4$, the formulas (\ref{mamuka1}) coincide with the formulas for the first two bending moments given in the second section. 
\end{rem}

The results on the complex invariant moments given in the preceding section enable us to obtain the following analog of Theorem \ref{main2} in the similar way as for the invariant bending moments. We omit a detailed proof since it is basically a repetition of the argument used to prove Theorem \ref{main2}.

\begin{thm} \label{main3}
	For $n\geq 3$, each invariant complex moment $J_{k,m}$ of poristic Steiner $n$-chains is algebraically dependent with the first two bending moments $I_1, I_2$. 
\end{thm}

We conclude this section by outlining an application of our results to the so-called feasibility problem for the radii of Steiner $4$-chains considered in \cite{kir} and \cite{bibdia}. Recall that this problem is formulated as follows. Given an ordered quadruple of positive numbers $(r_1, r_2, r_3, r_4)$ find out if there exists a Steiner $4$-chain with such radii of its circles in the given order. To the best of our knowledge, this natural problem was considered for the first time in \cite{kir}, where the author gave its solution in a concrete case but did not discuss the general case. A general algorithmic solution to this problem was given in \cite{bibdia}. We show that the verification of the criterion given in \cite{bibdia} can be simplified using Theorem \ref{thirdmomentrelation}. \\

Recall that the algorithm given in \cite{bibdia} involved several steps. First, one computes the first three moments $(I_1, I_2, I_3)$ of the quadruple $b_i = 1/r_i$ called the {\it actual moments of bends}. Then assuming the existence of a sought Steiner $4$-chain one uses the equations (\ref{moment41}) and (\ref{moment42}) to find the "virtual radii" $(\tilde R, \tilde r)$ of the "virtual parent circles". To this end it is sufficient to solve the system
$${2(a+A)= I_1, 3A^2 + 10Aa + 3a^2 = 2I_2},$$
which reduces to the following quadratic equation
\begin{equation} \label{parentradii}
	16a^2 - 8I_1a + (8I_2 - 3I_1^2) = 0.
\end{equation}
So the solutions of the above system can be expressed by explicit formulas, which yields the virtual parent radii $(\tilde R, \tilde r)$. If $\tilde R^2 - 6 \tilde R \tilde r + \tilde r^2 > 0$ then the pair $(\tilde R, \tilde r)$ is a feasible candidate for the radii of parent circles having the distance between their centers equal to $\tilde d = \sqrt{\tilde R^2 - 6 \tilde R \tilde r + \tilde r^2 }.$ Next, one has to verify that the given values of radii $r_i$ belong to the segment $[\tilde r_*, \tilde r^*]$. The final step of the algorithm given in \cite{bibdia,dia} required computing the "virtual bending moments" and comparing them with the the actual moments of bends $b_i = 1/r_i$. In view of our results this last step can be substituted by verifying that the actual moments of bends satisfy the relation (\ref{moment3byfirst2}), which is obviously simpler. In particular, if this relation is not satisfied then the given quadruple of positive numbers $r_i$ cannot be realized as the radii of any Steiner $4$-chain. We conclude this section by giving a simple example, where our algorithm easily yields a negative answer. \\

{\bf Example 2.}  If we are given a quadruple of radii $\rho = (1, 2, 3, 4)$, then the bends are $\varkappa = (1, 0.5, 0.33(3), 0.25).$ So $(I_1,I_2,I_3)=(2.083,1.4214,1.1766)$ and we see that (\ref{moment3byfirst2}) is not satisfied. So we conclude that the quadruple $(1, 2, 3, 4)$ is not realizable as the radii of a Steiner $4$-chain.

Analogous applications to the feasibility problem for the centers of Steiner chains can be given in terms of complex invariant moments but we do not dwell on these topics here for the sake of brevity.

\section{\textbf{Concluding remarks}}

First of all, it is natural to search for generalizations of our results to Steiner $n$-chains for arbitrary $n$. For small $n$, this seems feasible in the same way as above, using (\ref{yiuquadratics}) and the invariance of the first $n-1$ moments of curvatures \cite{scta}. For $n=6$, this is especially simple since one can find the invariant bending moments in a similar way by finding the bends in the axial symmetric $6$-chain as above. Specifically, using the symmetry of the axial chain we conclude that the collection of its bends is 

$$\frac{2}{R+d-r}, \frac{3R^2 + 3Rd - 2Rr - 3dr + 3r^2}{4Rr(R - r + d)}, \frac{3R^2 - 3Rd - 2Rr + 3dr + 3r^2}{4Rr(R - r - d)},$$ \\
$$\frac{2}{R-r-d}, \frac{3R^2 - 3Rd - 2Rr + 3dr + 3r^2}{4Rr(R - r - d)}, \frac{3R^2 + 3Rd - 2Rr - 3dr + 3r^2}{4Rr(R - r + d)}.$$ \\

Summing their powers one can obtain explicit formulas for the first five invariant bending moments. The resulting formulas are rather lengthy and will be published elsewhere.

Next, one can use the explicit formulas for the invariant bending moments to compute the canonical coordinates of the axial symmetric chain. As is easy to realize this reduces to computing the intersections of pairs of quadrics, which can be done explicitly. This in turn can be used to compute the complex moments for the symmetric axial chains.

Finally, one can search for the algebraic relations between complex invariants moments existence of which was shown above and use these relations for solving the feasibility problem for the centers of Steiner chains. \\

	\end{document}